\renewcommand{\Box}{\#}

\documentclass{article}
\newtheorem{theorem}{Theorem}[section]
\newtheorem{conjecture}[theorem]{Conjecture}
\newtheorem{corollary}[theorem]{Corollary}

\newtheorem{lemma}[theorem]{Lemma}

\def\cadre{$$\vcenter\bgroup\advance\hsize by -6em\noindent
             \refstepcounter{equation}\ignorespaces}
\makeatletter
\def\endcadre{\egroup\eqno(\theequation)$$\global\@ignoretrue}
\makeatother

\def\vertex{\circle*{0.4}}

\begin{document}

\begin{center}

{\LARGE {\bf A class of perfectly contractile graphs}}

\vspace{2em}

{\large{\bf Fr\'{e}d\'{e}ric Maffray\footnote{C.N.R.S.}, Nicolas
Trotignon\footnote{Supported by Universit\'e Pierre Mendes France,
Grenoble.}}}

{\small Laboratoire Leibniz-IMAG,  46 Avenue F\'{e}lix Viallet,\\
38031~Grenoble~Cedex, France.\\
(frederic.maffray@imag.fr, nicolas.trotignon@imag.fr)}

\end{center}
\begin{quote}
\small {\bf Abstract.} We consider the class ${\cal A}$ of graphs that
contain no odd hole, no antihole, and no ``prism'' (a graph consisting
of two disjoint triangles with three disjoint paths between them).  We
prove that every graph $G\in{\cal A}$ different from a clique has an
``even pair'' (two vertices that are not joined by a chordless path of
odd length), as conjectured by Everett and Reed [see the chapter
``Even pairs'' in the book {\it Perfect Graphs},
J.L.~Ram\'{\i}rez-Alfons\'{\i}n and B.A.~Reed, eds., Wiley
Interscience, 2001].  Our proof is a polynomial-time algorithm that
produces an even pair with the additional property that the
contraction of this pair yields a graph in ${\cal A}$.  This entails a
polynomial-time algorithm, based on successively contracting even
pairs, to color optimally every graph in ${\cal A}$.  This generalizes
several results concerning some classical families of perfect graphs.
\end{quote}

\section{Introduction}

A graph $G$ is \emph{perfect} if every induced subgraph $G'$ of $G$
satisfies $\chi(G')=\omega(G')$, where $\chi(G')$ is the chromatic
number of $G'$ and $\omega(G')$ is the maximum clique size in $G'$.
Berge {\cite{ber60,ber61,ber85}} introduced perfect graphs and
conjectured that {\it a graph is perfect if and only if it does not
contain as an induced subgraph an odd hole or an odd antihole of
length at least $5$}, where a \emph{hole} is a chordless cycle with at
least four vertices and an \emph{antihole} is the complement of a
hole.  We follow the tradition of calling \emph{Berge graph} any graph
that contains no odd hole and no odd antihole of length at least $5$.
This famous question (the Strong Perfect Graph Conjecture) was the
objet of much research (see the book \cite{ramree01}), until it was
finally proved by Chudnovsky, Robertson, Seymour and Thomas
\cite{CRST}: \emph{Every Berge graph is perfect}.  Moreover, a
polynomial-time algorithm was devised to decide if a graph is Berge
(hence perfect); it is due to Chudnovsky, Cornu\'ejols, Liu, Seymour
and Vu\v{s}kovi\'c \cite{CCLSV2002,CLV2002,CS2002}.

Despite those breakthroughs, some conjectures about Berge graphs
remain open.  An \emph{even pair} in a graph $G$ is a pair $\{x,y\}$
of non-adjacent vertices having the property that every chordless path
between them has even length (number of edges).  Given two vertices
$x,y$ in a graph $G$, the operation of \emph{contracting} them means
removing $x$ and $y$ and adding one vertex with edges to every vertex
of $G\setminus \{x,y\}$ that is adjacent in $G$ to at least one of
$x,y$; we denote by $G/xy$ the graph that results from this operation.
Fonlupt and Uhry \cite{fonuhr82} proved that \emph{if $G$ is a perfect
graph and $\{x,y\}$ is an even pair in $G$, then the graph $G/xy$ is
perfect and has the same chromatic number as $G$}.  In particular,
given a $\chi(G/xy)$-coloring $c$ of the vertices of $G/xy$, one can
easily obtain a $\chi(G)$-coloring of the vertices of $G$ as follows:
keep the color for every vertex different from $x,y$; assign to $x$
and $y$ the color assigned by $c$ to the contracted vertex.  This idea
could be the basis for a conceptually simple coloring algorithm for
Berge graphs: as long as the graph has an even pair, contract any such
pair; when there is no even pair find a coloring $c$ of the contracted
graph and, applying the procedure above repeatedly, derive from $c$ a
coloring of the original graph.  The algorithm for recognizing Berge
graphs mentioned at the end of the preceding paragraph can be used to
detect an even pair in a Berge graph $G$; indeed, it is easy to see
that two non-adjacent vertices $a,b$ form an even pair in $G$ if and
only if the graph obtained by adding a vertex adjacent only to $a$ and
$b$ is Berge.  Thus, given a Berge graph $G$, one can try to color its
vertices by keeping contracting even pairs until none can be found.
Then some questions arise: what are the Berge graphs with no even
pair?  What are, on the contrary, the graphs for which a sequence of
even-pair contractions leads to graphs that are trivially easy to
color?

\begin{figure}[htb]
\unitlength=0.3cm
\begin{center}\begin{tabular}{ccc}
\begin{picture}(6,13)
\multiput(0,0)(0,6){2}{\vertex}
\multiput(2,2)(0,2){2}{\vertex}
\multiput(4,0)(0,6){2}{\vertex}
\multiput(0,0)(0,6){2}{\line(1,0){4}}
\multiput(0,0)(4,0){2}{\line(0,1){6}}
\put(2,2){\line(0,1){2}}
\multiput(0,0)(2,4){2}{\line(1,1){2}}
\multiput(0,6)(2,-4){2}{\line(1,-1){2}}
\put(2,-1.5){\makebox(0,0){An odd prism: $\overline{C}_6$}}
\end{picture}
\quad\quad\quad & \quad\quad\quad
\begin{picture}(6,13)
\multiput(1,0)(4,0){2}{\vertex}
\multiput(0,2)(3,0){3}{\vertex}
\multiput(0,5)(3,0){3}{\vertex}
\multiput(1,7)(4,0){2}{\vertex}
\multiput(1,0)(0,7){2}{\line(1,0){4}}
\multiput(0,2)(3,0){3}{\line(0,1){3}}
\multiput(1,0)(2,5){2}{\line(1,1){2}}
\multiput(1,7)(2,-5){2}{\line(1,-1){2}}
\multiput(0,5)(5,-5){2}{\line(1,2){1}}
\multiput(0,2)(5,5){2}{\line(1,-2){1}}
\put(3,-1.5){\makebox(0,0){An odd prism}}
\end{picture}
\quad\quad\quad & \quad\quad\quad
\begin{picture}(6,13)
\multiput(0,0)(0,4){3}{\vertex}
\multiput(2,2)(0,2){3}{\vertex}
\multiput(4,0)(0,4){3}{\vertex}
\multiput(0,0)(0,8){2}{\line(1,0){4}}
\multiput(0,0)(4,0){2}{\line(0,1){8}}
\put(2,2){\line(0,1){4}}
\multiput(0,0)(2,6){2}{\line(1,1){2}}
\multiput(0,8)(2,-6){2}{\line(1,-1){2}}
\put(2,-1.5){\makebox(0,0){An even prism}}
\end{picture}
\end{tabular}\end{center}
\vspace{0.6cm}
\caption{Some prisms}\label{fig:prisms}
\end{figure}
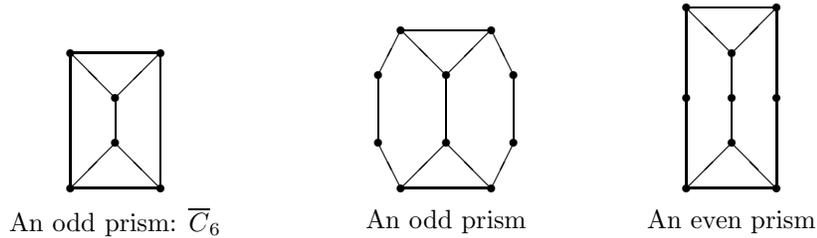

As a first step towards getting a better grasp on these questions,
Bertschi \cite{ber90} proposed the following definitions.  A graph $G$
is \emph{even-contractile} if either $G$ is a clique or there exists a
sequence $G_0, \ldots, G_k$ of graphs such that $G=G_0$, for $i=0,
\ldots, k-1$ the graph $G_i$ has an even pair $\{x_i, y_i\}$ such that
$G_{i+1}=G_i/x_iy_i$, and $G_k$ is a clique.  A graph $G$ is \emph{
perfectly contractile} if every induced subgraph of $G$ is
even-contractile.  This class is of interest because it turns out that
many families of graphs that are considered classical are perfectly
contractile (see Section~\ref{sec:more}).  Everett and Reed proposed a
conjecture aiming at a characterization of perfectly contractile
graphs.  To understand it, one more definition is needed: say that a
graph is a \emph{prism} if it consists of two vertex-disjoint
triangles (cliques of size $3$) with three vertex-disjoint paths
between them, and with no other edge than those in the two triangles
and in the three paths.  (Prisms were called stretchers in
\cite{linmafree97,linmaf01}).  Say that a prism is odd if these three
paths all have odd length.  For example the graph $\overline{C}_6$
that is the complement of a hole of length $6$ is an odd prism.  See
Figure {\ref{fig:prisms}}.  %

\begin{conjecture}[Everett and Reed \cite{epsbook,ree93}]
A graph is perfectly contractile if and only if it contains no odd
hole, no antihole of length at least $5$, and no odd prism.
\end{conjecture}
The if part of this conjecture remains open.  The only if part is not
hard to establish, but it requires some careful checking; this was
done formally in \cite{linmafree97}.  A weaker form was also proposed
by Everett and Reed:
\begin{conjecture}[Everett and Reed \cite{epsbook,ree93}]
If a graph contains no odd hole, no antihole of length at least $5$,
and no prism then it is perfectly contractile.
\end{conjecture}
We will prove here this second conjecture.  For this purpose, some
definitions and notation must be introduced.

Let $G=(V,E)$ be a graph.  Its complementary graph is denoted by
$\overline{G}$.  The subgraph induced by any $X\subseteq V$ is denoted
by $G[X]$.  We say that a vertex $u$ \emph{sees} a vertex $v$ when $u,
v$ are adjacent, else we say that $u$ \emph{misses} $v$.  For any $T
\subseteq V$, we let $N(T)$ denote the set of vertices of $G \setminus
T$ that see at least one vertex of $T$.  A vertex of $V\setminus T$ is
called \emph{$T$-complete} if it sees all vertices of $T$; then $C(T)$
denotes the set of $T$-complete vertices of $V\setminus T$.  We call
\emph{$T$-edge} any edge whose two vertices are $T$-complete.

A non-empty set $T\subseteq V$ is \emph{interesting} if
$\overline{G}[T]$ is connected (in short we will say that $T$ is
co-connected) and $G[C(T)]$ is not a clique (we view the empty set as
a clique, so $|C(T)|\ge 2$).  Note that a graph $G$ may fail to have
any interesting set; in that case it must be in particular that the
neighbourhood $N(v)$ of every vertex $v$ induces a clique ($v$ is
\emph{ simplicial}), for otherwise $\{v\}$ would be interesting; this
means that every connected component of $G$ is a clique, i.e., $G$ is
a disjoint union of cliques.

Paths in a graph may be described in two different ways: either we
list their vertices in order ($P=ab\cdots cd$); or, if $x,y$ are two
vertices of a given path $P$, we may use $P[x,y]$ or $P[y,x]$ to
denote the subpath of $P$ whose endvertices are $x,y$.  The length of
a path is its number of edges.  An edge between two vertices that are
not consecutive along the path is a \emph{chord}, and a path that has
no chord is \emph{chordless}.

\begin{figure}[htb]
\unitlength=0.26cm
\hspace{1cm}
\begin{picture}(18,5)
\multiput(2,0)(8,0){2}{\vertex}
\multiput(2,4)(8,0){2}{\vertex}
\multiput(0,2)(12,0){2}{\vertex}
\multiput(4,0)(2,0){3}{\vertex}
\multiput(4,4)(2,0){3}{\vertex}
\put(-1,1.5){$a$}
\put(2.5,3){$c$}
\put(9,3){$c'$}
\put(2.5,0.5){$d$}
\put(9,0.5){$d'$}
\put(12.5,1.5){$b$}
\multiput(2,0)(0,4){2}{\line(1,0){8}}
\multiput(2,0)(8,0){2}{\line(0,1){4}}
\multiput(0,2)(10,-2){2}{\line(1,1){2}}
\multiput(0,2)(10,2){2}{\line(1,-1){2}}
\end{picture}
\hfill
\begin{picture}(24,5)
\multiput(4,0)(8,0){2}{\vertex}
\multiput(4,4)(8,0){2}{\vertex}
\multiput(0,2)(2,0){2}{\vertex}
\multiput(14,2)(2,0){3}{\vertex}
\put(8,0){\vertex}
\multiput(5.6,4)(1.6,0){5}{\vertex}
\put(-1,1.5){$a$}
\put(1.5,2.5){$a'$}
\put(4.5,3){$c$}
\put(11,3){$c'$}
\put(4.5,0.5){$d$}
\put(11,0.5){$d'$}
\put(14,2.5){$b'$}
\put(18.5,1.5){$b$}
\multiput(4,0)(0,4){2}{\line(1,0){8}}
\multiput(4,0)(8,0){2}{\line(0,1){4}}
\multiput(2,2)(10,-2){2}{\line(1,1){2}}
\multiput(2,2)(10,2){2}{\line(1,-1){2}}
\put(0,2){\line(1,0){2}}
\put(14,2){\line(1,0){4}}
\end{picture}
\vspace{0.6cm}
\caption{Some $(a,b)$-snakes; the second one is
proper}\label{fig:snakes}
\end{figure}
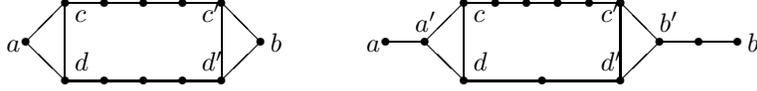

A \emph{snake} $S$ is a graph that consists of four disjoint chordless
paths $S_1= a\cdots a'$, $S_2=b\cdots b'$, $S_3=c\cdots c'$, $S_4=
d\cdots d'$, where $S_1, S_2$ may have length $0$ but $S_3, S_4$ have
length at least $1$, and such that the edge-set of $S$ is $E(S_1)\cup
E(S_2)\cup E(S_3) \cup E(S_4)\cup\{a'c, a'd, cd, b'c', b'd', c'd'\}$.
Note that $a'cd$ and $b'c'd'$ are triangles in $S$.  Vertices $a$ and
$b$ are the endvertices of the snake, and we may also say that $S$ is
an $(a, b)$-snake.  See Figure \ref{fig:snakes}.  A snake is
\emph{proper} if at least one of $S_1, S_2$ has length at least $1$.
An even pair $\{a,b\}$ in a graph $G$ is \emph{special} if the graph
$G$ contains no proper $(a, b)$-snake.

Let $\cal A$ be the class of graphs that contain no odd hole, no
antihole of length at least $5$, and no prism.  Our main result is the
following.  %
\begin{theorem}\label{thm:main}
In a graph $G\in{\cal A}$ that is not a disjoint union of cliques, let
$T$ be any interesting set.  Then $C(T)$ contains a special even pair
of $G$.
\end{theorem}
The proof of this theorem will be given in Section~\ref{sec:proof}.
Some technical lemmas will be given in Section~\ref{sec:techlems}.

Special even pairs are valuable because of the following result (see
\cite{epsbook}):
\begin{theorem}[\cite{epsbook}]\label{thm:sep}
Let $G$ be a graph and $\{x,y\}$ be an even pair of $G$.  Then:
\begin{enumerate}
\item
If $G$ contains no odd hole, then $G/xy$ contains no odd hole;
\item
If $G$ contains no antihole of length at least $5$, then $G/xy$
contains no antihole of length at least $5$ except possibly
$\overline{C}_6$;
\item
If $G$ contains no prism and no proper $(x,y)$-snake, then $G/xy$
contains no prism.
\end{enumerate}
\end{theorem}
Note that $\overline{C}_6$ (the complement of the $6$-vertex hole) is
itself a prism.
\begin{corollary}\label{cor:main}
Every graph $G\in{\cal A}$ either is a clique or has a special even
pair.  Every graph $G\in{\cal A}$ is perfectly contractile.
\end{corollary}
\emph{Proof.} If $G$ is a disjoint union of cliques and not itself a
clique, then any two non-adjacent vertices $x,y$ form a special even
pair.  If $G$ is not a disjoint union of cliques, apply Theorem
\ref{thm:main} to get a special even pair $\{x,y\}$.  In either case,
Theorem~\ref{thm:sep} implies that $G/xy$ is in ${\cal A}$, which by
induction entails the second sentence of the corollary.  \hfill$\Box$

\section{Some technical lemmas}\label{sec:techlems}

First we recall a nice lemma due to Roussel and Rubio \cite{rourub01},
also proved independently by Chudnovsky, Robertson, Seymour, Thomas,
and Thomassen (but not published) and used intensively in \cite{CRST}.
We propose here our own proof of that lemma, which we believe is
simpler and shorter than Roussel and Rubio's.  Then we derive some
lemmas of a similar flavor.

For any chordless path $P=xx'\cdots y'y$ of length at least $2$ in
$G$, let $P^*$ denote the interior of $P$, i.e., the path $P\setminus
\{x, y\}$.  Following \cite{CRST}, we say that a pair $\{u, v\}$ of
non-adjacent vertices of $V \setminus P$ is a \emph{leap} for $P$ if
$N(u)\cap P=\{x, x', y\}$ and $N(v)\cap P=\{x, y', y\}$.  Note that in
that case $P^*\cup\{u,v\}$ is a chordless path of the same length as
$P$.

\begin{lemma}[\cite{rourub01}]\label{lem:w}
In a Berge graph $G=(V,E)$, let $P,T\subset V$ be disjoint sets such
that $P$ induces a chordless path, $T$ induces a co-connected
subgraph, and the endvertices of $P$ are $T$-complete.  Then one of
the following four outcomes holds:
\begin{enumerate}
\setcounter{enumi}{-1}
\item \label{oeven}
$P$ has even length  and has an even number of $T$-edges;
\item \label{oodd}
$P$ has  odd length and has an odd number of $T$-edges;
\item \label{oleap}
$P$ has odd length at least $3$ and there is a leap for $P$ in $T$;
\item \label{ohop}
$P$ has length $3$ and its two internal vertices are the endvertices
of a chordless odd path of $\overline{G}$ whose interior is in $T$.
\end{enumerate}
\end{lemma}
\emph{Proof.} We prove the lemma by induction on $|P\cup T|$.  If $P$
has length $0$ or $1$ then we have outcome~\ref{oeven} or~\ref{oodd}.
So let us assume that $P$ has length at least $2$.  Put $P=xx'\cdots
y'y$.  We distinguish between two cases.

\emph{Case 1: There is no $T$-complete vertex in $P^*$.} If $P$ has
even length, we have outcome~\ref{oeven}.  So we may assume that $P$
has odd length.  If $|T|=1$, then $T\cup P$ induces an odd hole.  So
$|T| \ge 2$.  Let us suppose that outcomes~\ref{oleap} and~\ref{ohop}
do not hold for $P$.  Therefore, and by induction, we know that for
every co-connected proper subset $U$ of $T$ there is an odd number of
$U$-edges in $P$.  Note that, for any $t\in T$, any $T\setminus
\{t\}$-complete vertex of $P^*$ misses $t$.

\emph{Case 1.1: $T$ induces a stable set.} Let $t$ be any vertex of
$T$.  We claim that $N(t)\cap P = \{x, x', y\}$ or $N(t)\cap P =
\{x,y',y\}$.  Call \emph{$t$-segment} of $P$ any subpath of $P$ whose
endvertices see $t$ and whose internal vertices miss $t$.  Since $x,
y$ see $t$, $P$ is edge-wise partitioned into its $t$-segments.  Also,
we know that there is an odd number of $T \setminus \{t\}$-edges in
$P$.  So there is a $t$-segment $P[r, s]$ of $P$ that has an odd
number of $T \setminus \{t\}$-edges.  Assume that $x, r, s, y$ appear
in this order along $P$.  Note that $P[r, s]$ has length at least $2$,
for otherwise one of $r,s$ would be a $T$-complete vertex in $P^*$.
So $P[r,s] \cup\{t\}$ induces an even hole and $P[r,s]$ has even
length.  Let $r', s'$ respectively be the $T \setminus \{t\}$-complete
vertices of $P[r,s]$ closest to $r$ and to $s$.  Note that $P[r',s']$
contains all the $T \setminus \{t\}$-edges of $P[r, s]$, so the
induction hypothesis, applied to $P[r', s']$ and $T\setminus\{t\}$,
implies that $P[r',s']$ has odd length.  Thus exactly one of the paths
$P[r,r']$ and $P[s',s]$ has odd length.  Assume that $P[r, r']$ has
odd length.  So $r\neq r'$ and so $x \neq r$.

If $r'$ misses $y$, then the odd path $P_1=P[r',r]\cup \{t,y\}$ is
chordless.  The endvertices $r', y$ of $P_1$ are $T \setminus
\{t\}$-complete, and $P_1$ has no $T \setminus \{t\}$-edge since
$P_1^*$ contains no $T\setminus\{t\}$-complete vertex.  Thus, by the
induction hypothesis, we must have outcome~\ref{oleap} or~\ref{ohop}
for $P_1$ and $T \setminus \{t\}$; but these are impossible since $t$
has no neighbour in $T \setminus \{t\}$.  So $r'$ must see $y$,
meaning that $s'=s=y$ and $r'= y'$.

If $r$ misses $x$, then the odd path $P_2=P[r',r]\cup \{t,x\}$ is
chordless.  The endvertices $r', x$ of $P_2$ are $T \setminus
\{t\}$-complete, and $P_2$ has no $T \setminus \{t\}$-edge since
$P_2^*$ contains no $T\setminus\{t\}$-complete vertex.  Thus, by the
induction hypothesis, we must have outcome~\ref{oleap} or~\ref{ohop}
for $P_2$ and $T \setminus \{t\}$; but these again are impossible
since $t$ has no neighbour in $T \setminus \{t\}$.  So, $r$ must see
$x$, meaning that $r = x'$.  Thus we have $N(t) \cap P = \{x, x',
y\}$.  Likewise if $P[s', s]$ has odd length then $N(t) \cap P = \{x,
y', y\}$.  So the above claim is proved, for every $t\in T$.

Now, since $y'$ is not $T$-complete there is a vertex $u\in T$ such
that $N(u) \cap P = \{x,x',y\}$, and since $x'$ is not $T$-complete
there is a vertex $v\in T$ such that $N(v) \cap P = \{x,y',y\}$.
Thus, $\{u, v\}$ is a leap in $T$ for $P$.

\emph{Case 1.2: $T$ does not induce a stable set.} Let $Q=u\cdots v$ be
a longest path of $\overline{G}[T]$.  So $Q$ has length at least $2$
(since $T$ is not a stable set), and $T\setminus \{u\}$ and $T
\setminus \{v\}$ are co-connected sets.  So we know that $P$ has an
odd number of $T\setminus \{u\}$-edges and an odd number of $T
\setminus \{v\}$-edges.  Note that a $T \setminus \{u\}$-edge and a $T
\setminus \{v\}$-edge have no common vertex, for otherwise there would
be a $T$-complete vertex in $P^*$.  In particular all $T \setminus
\{u\}$-edges and $T\setminus\{v\}$-edges are different.

Suppose that $Q$ has even length.  Let $x_u x'_u$ be a $T\setminus
\{u\}$-edge of $P$ and $y'_v y_v$ be a $T\setminus \{v\}$-edge of $P$
such that, without loss of generality, $x, x_u, x'_u, y'_v, y_v, y$
appear in this order on $P$.  If $x'_u$ misses $y'_v$ then $\{x'_u,
y'_v\} \cup Q$ induces an odd antihole.  If $x\neq x_u$ then $\{x_u,
y'_v\} \cup Q$ induces an odd antihole.  If $y_v \neq y$ then $\{x'_u,
y_v\} \cup Q$ induces an odd antihole.  It follows that $P = x_u x'_u
y'_v y_v$, but then $P\cup Q$ induces an odd antihole.  Thus $Q$ has
odd length (at least $3$).

Suppose that $T\setminus\{u,v\}$ is not co-connected.  Since $T 
\setminus \{u\}$ and $T\setminus\{v\}$ are co-connected, there exists
a vertex $w$ in a connected component of $\overline{G}[T] \setminus
\{u,v\}$ that does not contain $Q^*$ and such that $w$ sees in
$\overline{G}$ at least one of $u,v$; but then $Q\cup\{w\}$ induces in
$\overline{G}[T]$ either a chordless path longer than $Q$ or an odd
hole, a contradiction.  So $T\setminus\{u,v\}$ is co-connected.

Now we know that there is an odd number of $T\setminus \{u, v\}$-edges
in $P$.  Recall that $P$ has an odd number of $T\setminus\{u\}$-edges,
an odd number of $T\setminus\{v\}$-edges, and that these are
different, so these account for an even number of $T\setminus \{u,
v\}$-edges; thus $P$ has at least one $T\setminus \{u, v\}$-edge $x''
y''$ that is neither a $T\setminus \{u\}$-edge nor a $T \setminus
\{v\}$-edge.  We may assume that $x,x'',y'',y$ appear in this order
along $P$ and that $y''\in P^*$.  So $y''$ misses one of $u, v$, say
$v$.  Then $y''$ sees $u$, for otherwise $Q \cup \{y''\}$ would induce
an odd antihole.  Then $x''$ misses $u$, for otherwise $x''y''$ would
be a $T\setminus \{v\}$-edge.  Then $x''$ sees $v$, for otherwise
$Q\cup \{x''\}$ would induce an odd antihole.  Then $x''=x'$ for
otherwise $Q\cup\{x'', y'', x\}$ would induce an odd antihole, and
similarly $y''=y'$.  So $P=xx''y''y$ and $Q\cup\{x'', y''\}$ is a
chordless odd path of $\overline{G}$, and we have outcome~\ref{ohop}.

\emph{Case 2: There is a $T$-complete vertex in $P^*$.} Let $z$ be such
a vertex.  By induction, we can apply the lemma to the path $P_{xz}=
P[x\cdots z]$ and $T$.  If $P_{xz}$ is odd and there is a leap $\{u,
v\}$ for $P_{xz}$ in $T$, then $P_{xz}^* \cup \{u, v, y\}$ induces an
odd hole.  If $P_{xz}$ has length $3$ and its two internal vertices
are the endvertices of a chordless odd path $M$ of $\overline{G}$
whose interior is in $T$, then $M\cup\{y\}$ induces an odd antihole.
So it must be that the number of $T$-edges in $P_{xz}$ and the length
of $P_{xz}$ have the same parity.  The same holds for $P[z \cdots y]$.
Then the number of $T$-edges in $P$ and the length of $P$ have the
same parity, and we have outcome~\ref{oeven} or~\ref{oodd}.
\hfill$\Box$

Let ${\cal A}'$ be the class of graphs that contain no odd hole, no
antihole of length at least $5$, and no odd prism.  Clearly ${\cal
A}\subset{\cal A}'$.  %
\begin{lemma}\label{lem:wp}
In a graph $G=(V,E)\in{\cal A}'$, let $P,T\subset V$ be disjoint sets
such that $P$ induces a chordless path, $T$ induces a co-connected
subgraph, and the endvertices of $P$ are in $C(T)$.  Then the number
of $T$-edges in $P$ has the same parity as the length of $P$.  In
particular if $P$ has odd length then some internal vertex of $P$ is
in $C(T)$.
\end{lemma}
\emph{Proof.} Apply Lemma~\ref{lem:w}.  Observe that in
outcome~\ref{oleap} of Lemma~\ref{lem:w}, the set $P\cup\{u,v\}$
induces an odd prism, and that in outcome~\ref{ohop}, letting $M$
denote a chordless odd path of $\overline{G}$ whose interior is in $T$
and whose endvertices are the two internal vertices of $P$, $P\cup M$
induces an antihole of length at least $6$ in $G$.  Thus we must have
outcome~\ref{oeven} or~\ref{oodd}.  \hfill$\Box$

\begin{lemma}\label{lem:wh}
In a graph $G=(V,E)\in{\cal A}'$, let $H,T\subset V$ be disjoint sets
such that $H$ induces a hole, $T$ induces a co-connected subgraph, and
at least two non-adjacent vertices of $H$ are in $C(T)$.  Then the
number of $T$-edges in $H$ is even.
\end{lemma}
\emph{Proof.} Let $x, y$ be two non-adjacent vertices of $H\cap C(T)$.
Then the two $(x,y)$-paths contained in $H$ have the same parity, so
Lemma~\ref{lem:wp} implies that the numbers of $T$-edges they contain
have the same parity.  Thus $H$ has an even number of $T$-edges.
\hfill$\Box$

In a graph $G$, we say that three paths $P_1$, $P_2$, $P_3$ induce a
$\Delta P(x_1x_2x_3, x)$ (where $x, x_1, x_2, x_3$ are four distinct
vertices) if each $P_i$ is a chordless $(x_i, x)$-path, $x_1, x_2,
x_3$ induce a triangle, at least two of the paths have length at least
$2$, and the $\Delta P$ has no other edge than those in the three
paths and in the triangle.  If a graph $G$ contains such a
configuration then two of the three paths have the same parity, and
so, since at least one of them has length at least $2$, the union of
these two paths induces an odd hole.  Thus: %
\begin{lemma}\label{lem:dp}
In a graph $G$ that contains no odd hole, there is no $\Delta P$
configuration.  \hfill$\Box$
\end{lemma}

\begin{lemma}\label{lem:sgp}
In a graph $G=(V,E)\in {\cal A}$, let $S,T\subset V$ be disjoint sets
such that $S$ induces an $(a,b)$-snake, $T$ induces a co-connected
subgraph, and $a,b\in C(T)$.  Then for each triangle of $S$ at least
two vertices of the triangle are in $C(T)$.
\end{lemma}
\emph{Proof.} We use the notation for snakes given above.  We first
claim that every $t\in T$ sees at least two vertices of the triangle
$a'cd$.  Let $S_1^t$ be a chordless $(a', t)$-path contained in
$S_1\cup \{t\}$, and $S_2^t$ be a chordless $(b', t)$-path contained
in $S_2\cup \{t\}$.  Put $H=S_3\cup S_4$.

If $t$ has no neighbour in $H$ then $H\cup S_1^t\cup S_2^t$ induces a
prism (which has the same two triangles as $S$), a contradiction.

Suppose that $t$ has exactly one neighbour $h$ in $H$.  If $h\notin
\{c,d\}$, then $H \cup S_1^t$ induces a $\Delta P(a' c d, h)$, and if
$h\in\{c,d\}$, then $H \cup S_2^t$ induces a $\Delta P(b' c' d', h)$,
in either case a contradiction to Lemma~\ref{lem:dp}.

Suppose that $t$ has exactly two neighbours $h, h'$ in $H$ and these
are adjacent.  Call $c''$ (resp.~$d''$) the neighbour of $c$ along
$S_3$ (resp.~of $d$ along $S_4$).  If the pair $\{h,h'\}$ is none of
the two pairs $\{c,c''\}, \{d,d''\}$ then $H\cup S_1^t$ induces a
prism (with triangles $a' c d$, $t h h'$).  If $\{h, h'\}=\{c,
c''\}$, then either $S_1^t\cup \{c\}\cup S_4\cup S_2^t$ induces a
$\Delta P(a'cd, t)$, a contradiction, or $t$ sees $a'$ (thus $t$ sees
two of $a', c,d $ as desired).  If $\{h, h'\}=\{d, d''\}$, similarly
$t$ sees $a'$ and $d$.

Finally, suppose that $t$ sees two non-adjacent vertices of $H$.  Let
$h$ and $h'$ respectively be the vertices of $H\cap N(t)$ closest to
$c$ along $H\setminus\{d\}$ (resp.~closest to $d$ along $H\setminus
\{c\}$).  Call $H[c, h]$ the chordless $(c, h)$-path in $H \setminus
\{d\}$, and call $H[d, h']$ the chordless $(d, h')$-path in $H
\setminus \{c\}$.  Then, since $S_1^t\cup H[c, h] \cup H[d, h']$
cannot induce a $\Delta P(a' c d, t)$, it must be that at least two of
the three paths $S_1^t, H[c, h]\cup\{t\}, H[d, h']\cup\{t\}$ have
length $1$, so $t$ sees at least two of $a', c, d$.

Suppose now that the lemma fails for the first triangle of the snake,
that is, there exist vertices $\alpha,\beta\in \{a', c, d\}$ and
vertices $x,y\in T$ such that $x$ misses $\alpha$ (and thus sees
$\beta$) and $y$ misses $\beta$ (and thus sees $\alpha$).  Since $T$
is co-connected, there is a chordless $(x, y)$-path $R$ in
$\overline{G}(T)$, and we can choose $x,y$ so that $R$ is as short as
possible; it follows that all the internal vertices of $R$ see both
$\alpha, \beta$ in $G$.  But now $R\cup \{\alpha, b, \beta\}$ induces
an antihole of length at least $5$ in $G$, a contradiction.  Thus at
least two vertices among $a', c, d$ are in $C(T)$.  The same holds for
the other triangle of $S$.  \hfill$\Box$

\begin{lemma}\label{lem:sgt}
In a graph $G=(V,E)\in {\cal A}$, let $H, P, T\subset V$ be pairwise
disjoint sets such that $H$ induces an even hole, $P$ induces a
chordless $(x,y)$-path, $T$ induces a co-connected subgraph, and there
are two disjoint edges $ab, cd$ of $H$ such that the set of edges
between $H$ and $P$ is $\{ax, bx\}$ and $c,d,y$ are in $C(T)$.  Then
at least one of $a,b$ is in $C(T)$.
\end{lemma}
\emph{Proof.} Assume that $a,c,d,b$ lie in this order along $H$, and
call $P_1$ the chordless $(a,c)$-path contained in $H\setminus
\{b,d\}$ and $P_2$ the chordless $(b,d)$-path contained in $H
\setminus\{a,c\}$.

Suppose that the lemma does not hold: there exists a vertex $u\in
T\setminus N(a)$ and a vertex $v\in T\setminus N(b)$.  Let $a_u$ be
the vertex of $P_1\cap N(u)$ closest to $a$, let $b_u$ be the vertex
of $P_2\cap N(u)$ closest to $b$, and let $x_u$ be the vertex of
$P\cap N(u)$ closest to $x$.  If $a_u=c$ and $b_u=d$ then $P_1\cup
P_2\cup P[x,x_u]\cup\{u\}$ induce a prism (with triangles $abx$
and $ucd$), a contradiction.  If either $a_u\not=c$ or $b_u\not=d$,
then the three paths $P_1[a, a_u], P_2[b, b_u], P[x, x_u]$ have no
edge between them (other than the edges of the triangle $abx$), so, by
Lemma~\ref{lem:dp}, at least two of them have length $0$; this means
that $u$ sees $b$ and $x$.  Similarly $v$ sees $a$ and $x$.

Since $T$ is co-connected, there exists a chordless $(u,v)$-path $R$
in $\overline{G}(T)$.  We choose $u,v$ that minimize the length of
$R$, so the internal vertices of $R$ (if any) see both $a,b$.  If some
vertex $w\in\{y, c, d\}$ misses both $a,b$ then $R\cup\{a,b,w\}$
induces an antihole of length at least $5$, a contradiction.  In the
remaining case we must have $y=x$, $ac\in E, bd\in E$; but then
$R\cup\{a,b,c,d,y\}$ induces an antihole of length at least $7$.
\hfill$\Box$

\begin{lemma}\label{lem:dpg}
In a graph $G=(V,E)\in {\cal A}$, let $H, P, T\subset V$ be such that
$H$ induces a hole, $P$ induces a chordless $(x,y)$-path, $T$ induces
a co-connected subgraph disjoint from $H\cup P$, $H\cup P$ is
connected, and there are adjacent vertices $u,v\in H$ such that $x\neq
u,v$ and $u, v, x\in C(T)$.  Then, either some vertex of $P$ sees at
least one of $u, v$, or some vertex of $H\setminus\{u,v,x\}$ is in
$C(T)$.
\end{lemma}
\emph{Proof.} Note that $H$ and $P$ are not assumed to be disjoint.  We
may even have $P \subset H$.

Suppose that the lemma does not hold, and consider a counterexample
with $|H\cup P|$ minimal.

If $x\in H$, then, by Lemma~\ref{lem:wh}, $H$ has an even number of
$T$-edges; thus $H$ has a $T$-edge different from $uv$, so there is a
vertex of $H\setminus \{u,v,x\}$ in $C(T)$, and the triple $(H, P, T)$
is not a counterexample to the lemma.  Therefore we may assume
$x\notin H$.

Let $x'$ be the vertex of $P$ that has a neighbour in $H$ and is
closest to $x$ along $P$.  So $P[x,x'] \cap H= \emptyset$.  Let $u'$
(resp.~$v'$) be the vertex of $H\cap N(x')$ closest to $u$ along
$H\setminus \{v\}$ (resp.~to $v$ along $H\setminus \{u\}$).  By the
assumption, $u'\not=u$ and $v'\not= v$.  Call $H[u, u']$ the path from
$u$ to $u'$ in $H\setminus\{v\}$, and call $H[v, v']$ the path from
$v$ to $v'$ in $H\setminus\{u\}$.

Suppose $u'= v'$.  The paths $P[x, x']\cup H[u', u]$ and $P[x, x']\cup
H[v', v]$ are chordless and one of them is odd.  By
Lemma~\ref{lem:wp}, this odd path must have an odd number of
$T$-edges.  Hence, there is at least one vertex $x''$ of $C(T)$ in
$P[x, x']\cup H\setminus\{u,v,x\}$.  If $x''$ is in $H$, the triple
$(H, P, T)$ is not a counterexample.  So $x''$ is in $P[x, x']
\setminus \{x\}$; but then the hole $H$, the path $P[x'', x']$, and
the set $T$ form a counterexample with $|H\cup P[x'', x']|<|H\cup P|$,
a contradiction.

Suppose $u'\neq v'$ and $u'v'\in E$.  Then we can apply
Lemma~\ref{lem:sgt} to the hole $H$, the path $P[x, x']$, and the set
$T$: so at least one of $u', v'$ is in $C(T)$, a contradiction.

Suppose $u' \neq v'$ and $u'v' \notin E$.  Consider the hole $H'$
induced by $H[u,u']\cup H[v,v']\cup\{x'\}$.  Then $H', P[x,x'], T$
form a counterexample to the lemma with $|H'\cup P[x,x']| < |H\cup
P|$, a contradiction.  \hfill$\Box$

Recall that a graph $G$ is \emph{weakly triangulated} \cite{hay85} if
$G$ and $\overline{G}$ contain no hole of length at least~$5$.
\begin{lemma}\label{lem:wt}
In  a  weakly triangulated  graph  $G=(V,E)$,  let  $P,T\subset V$  be
disjoint sets such that $P$ induces a chordless $(x,y)$-path of length
at least  $3$, $T$ induces a  co-connected subgraph, and  $x,y$ are in
$C(T)$.  Then at least one internal vertex of $P$ is in $C(T)$.
\end{lemma}
\emph{Proof.} Observe that no vertex $t\in T$ misses two consecutive
vertices $u,v$ of $P$, for otherwise $P\cup\{t\}$ contains a hole, of
length at least $5$, containing $u,v,t$.  Now let $v$ be an internal
vertex of $P$ that sees the most vertices of $T$.  If $v\in C(T)$ we
are done, so assume that there is a vertex $t\in T\setminus N(v)$.
Call $u,w$ the neighbours of $v$ on $P$, with $u\in P[x, v]$, $w\in
P[v, y]$.  Then both $tu, tw$ are edges, by the observation.  We may
assume that $u\not=x$ (else, by symmetry, $w\not=y$).  By the choice
of $v$, and since $t$ sees $u$ and misses $v$, there is a vertex
$t'\in T$ that sees $v$ and misses $u$.  Since $T$ is co-connected,
there is a chordless $(t,t')$-path $R$ in $\overline{G}[T]$, and we
choose $t,t'$ so that $R$ is as short as possible; so, and by the
observation, the internal vertices of $R$ see both $u,v$.  If $u$
misses $x$ then $R\cup \{u, v, x\}$ induces an antihole of length at
least $5$.  So $u$ sees $x$.  Likewise, $v$ sees $y$, for otherwise
$R\cup\{u,v,y\}$ induces an antihole of length at least $5$.  But then
$R\cup\{x, u, v, y\}$ induces an antihole of length at least $6$, a
contradiction.  \hfill$\Box$

\section{Proof of Theorem \protect{\ref{thm:main}}}\label{sec:proof}

We now prove Theorem~\ref{thm:main}.  Let $G=(V,E)$ be a graph that
contains no odd hole, no antihole of length at least $5$, and no prism
and that is not a disjoint union of cliques.  We proceed by induction
on $|V|$.  The smallest possible value is $|V|=3$, in which case $G$
is a three-vertex path and the desired conclusion is obvious.  Now let
us assume $|V|\ge 4$.  Along this proof we will also make several
remarks, for further reference, concerning the complexity of finding
various sets and paths relevant to the proof.

First observe that we need only prove the theorem for (inclusion-wise)
\emph{maximal} interesting sets.  For if $T$ is an interesting set in
$G$, then for any maximal interesting set $T'$ with $T\subset T'$ we
have $C(T')\subseteq C(T)$, thus any special even pair of $G$ in
$C(T')$ is also in $C(T)$.

So let $T$ be a maximal interesting set.  We observe that for every
vertex $z\in V\setminus (T\cup C(T))$, \emph{the set $N(z)\cap C(T)$
induces a clique in $G$}, for otherwise $T\cup\{z\}$ would be a larger
interesting set because $T\cup\{z\}$ is co-connected (since $T$ is
co-connected and $z\not\in C(T)$) and $C(T\cup \{z\})=C(T)\cap N(z)$.

{\bf Remark 1.} One can determine a maximal interesting set in
polynomial time.  Start from any non-simplicial vertex $v$ and put
$T:= \{v\}$.  (It may take time $O(|V||E|)$ to find such a $v$.)  As
long as there exists a vertex $w\in V\setminus (T\cup C(T))$ such that
$N(w)\cap C(T)$ does not induce a clique in $G$, put $T:=T\cup\{w\}$
and iterate.  At termination $T$ is a maximal interesting set.  There
may be $O(|V|)$ iterations, each taking time $O(|V||E|)$, so the time
to find such a set $T$ is $O(|V|^2|E|)$.

Let us call \emph{outer path} any chordless path whose endvertices
are non-adjacent vertices of $C(T)$ and whose internal vertices are in
$V\setminus (T\cup C(T))$.

Observe that \emph{there is no outer path of odd length}, by Lemma
\ref{lem:wp}.  Moreover, if $P$ is any outer path of even length, then
its length is at least $4$, for if it was $2$ the internal vertex $z$
of $P$ would be such that $N(z)\cap C(T)$ is not a clique, since it
contains the endvertices of $P$, a contradiction to the maximality of
$T$.  We now distinguish between two cases.

{\bf Case 1:} \emph{There is no outer path at all.}

Let $\{a,b\}$ be any special even pair of the graph $G[C(T)]$.  Recall
that $G[C(T)]$ is not a clique; so, such a pair exists, by the
induction hypothesis if $G[C(T)]$ is not a disjoint union of cliques,
trivially if $G[C(T)]$ is a disjoint union of cliques.  Consider any
chordless $(a, b)$-path $P$ in $G$.  If $P$ has a vertex $t\in T$ then
$P=atb$, so it has length $2$.  If $P\cap T=\emptyset$, it must be
that $P$ lies entirely in $C(T)$, for otherwise $P$ would contain an
outer path; so $P$ has even length.  Therefore $\{a, b\}$ is an even
pair of $G$.  Moreover, if there exists a proper $(a, b)$-snake $S$ in
$G$ (with the above notation for snakes), then the two $(a, b)$-paths
$P_1=S_1\cup S_3\cup S_2$ and $P_2=S_1\cup S_4\cup S_2$ have length at
least $4$, and so $P_1, P_2$ lie entirely in $C(T)$ by the same
argument as precedently with $P$; so $S$ lies entirely in $C(T)$, a
contradiction.  It follows that $\{a, b\}$ is a special even pair of
$G$.

{\bf Case 2:} \emph{There exists an (even) outer path.}

Let $\alpha z_1 \cdots z_n \beta$ be a \emph{shortest} outer path.
Its length is $n+1$.  Note that $n$ is odd and that $n\ge 3$ as
pointed out above.  Put $Z=z_1\cdots z_n$.  Define:
\begin{eqnarray*}
A&=&\{v\in C(T)\mid vz_1\in E, vz_i\not\in E \ (i=2,\ldots, n)\},\\
B&=&\{v\in C(T)\mid vz_n\in E, vz_i\not\in E \ (i=1,\ldots, n-1)\}.
\end{eqnarray*}
Note that $A$ is not empty, because $\alpha\in A$, and that $A$
induces a clique, because $A\subseteq N(z_1)\cap C(T)$.  Likewise $B$
is a non-empty clique.  Clearly, $A\cap B=\emptyset$.  Moreover, there
is no edge $uv$ with $u\in A, v\in B$, for otherwise $u, z_1, \ldots,
z_n, v$ would induce an odd hole.

{\bf Remark 2.} Finding a shortest outer path (or concluding that
there is no outer path at all) and, if there is any, finding the
corresponding sets $A,B$ can be done in in polynomial time.  Indeed it
suffices, for every pair of non-adjacent vertices $u, v\in C(T)$, to
look for a shortest $(u, v)$-path in $G\setminus(T\cup C(T) \setminus
\{u, v\})$, and to take the shortest of them all (if any) over all
pairs $u,v$.  Looking for a shortest path takes time $O(|E|)$ for each
pair $u,v$.  In total, since $T$ may have size $O(|V|)$, finding a
shortest outer path or concluding that there is none may take time
$O(|V|^2|E|)$.

We now show that some well-chosen vertices $a\in A$, $b\in B$ form a
special even pair of $G$.  This is established in Lemmas
\ref{lem:abi}--\ref{lem:absep}.  The outline of the proof from here on
is quite similar to that in \cite{linmaf01}.

\begin{lemma}\label{lem:abi}
$C(T)\cap N(Z)\subseteq A\cup B\cup C(T\cup A\cup B)$.
\end{lemma}
\emph{Proof.} Pick any $w\in C(T)\cap N(Z)$; so there exists an edge
$z_iw$ with $z_i\in Z$ ($1\le i\le n$).  If $w\in C(T\cup A\cup B)$ we
are done, so let us assume that there is a vertex $u\in A\cup B$ with
$uw\not\in E$, say $u\in A$.  Let $i$ be the smallest integer with
$z_iw\in E$.  Then $uz_1\cdots z_i w$ is an outer path, of length
$i+1$, so we must have $i=n$, and so $w\in B$.  \hfill$\Box$

\begin{lemma}\label{lem:pq}
Consider any odd chordless $(u,v)$-path $P= uu'\cdots v'v$, with $u\in
A$ and $v\in B$.  Then exactly one of $u'\in A$ or $v'\in B$ holds.
\end{lemma}
\emph{Proof.} Note that $P$ has length at least $3$, since there is no
edge between $A$ and $B$; also $P$ clearly contains no vertex of $T$
and no vertex of $C(T\cup A\cup B)$.

Suppoe that neither $u'\in A$ nor $v'\in B$ holds.  We will show that
this leads to a contradiction.  We claim that:
\begin{eqnarray}
\mbox{The only edges from $Z$ to $P\cap C(T)$ are $z_1u$ and
$z_nv$.}\label{eq:eqct}
\end{eqnarray}
For if $zw$ is an edge with $z\in Z$ and $w\in P\cap C(T)$, then,
since $w\notin C(T\cup A\cup B)$ as observed above, and by Lemma
\ref{lem:abi}, we have $w\in A\cup B$.  Since $A$ is a clique, the
case $w\in A$ means either $w=u$ (so $zw=z_1u$ as claimed) or $w=u'$
(so $u'\in A$, which we have excluded).  The case $w\in B$ is similar.
Thus (\ref{eq:eqct}) holds.

Let us mark those vertices of $P$ that have a neighbour in $Z$ (in
particular the vertices of $P\cap Z$ are marked).  Call \emph{
$Z$-segment} of $P$ any subpath of $P$, of length at least $1$,
whose endvertices are marked and whose internal vertices are not
marked.  By (\ref{eq:eqct}) the marked vertices of $P\setminus\{u,v\}$
are all in $P\setminus C(T)$.  Since $u,v$ are marked, $P$ is
(edge-wise) partitioned into its $Z$-segments.  Also some internal
vertex of $P$ must be marked, for otherwise $Z\cap P=\emptyset$ and
$Z\cup P$ induces an odd hole; so \emph{$P$ has at least two
$Z$-segments.}

By Lemma~\ref{lem:wp}, we know that \emph{$P$ has an odd number of
$T$-edges.}

It follows from the conclusion of the preceding two paragraphs that
there exists a $Z$-segment $Q$ of $P$ that contains an odd number of
$T$-edges, and that $Q$ does not contain both $u,v$, say $Q$ does not
contain $v$.  Call $w, x$ the endvertices of $Q$, and call $w', x'$
respectively the vertices of $Q\cap C(T)$ closest to $w$ and to $x$
along $Q$, so that $u, w, w', x', x, v$ lie in that order along $P$.
We have $w' \not=x'$ since $Q$ contains at least one $T$-edge.  For
the sake of clarity, note that $Q\cap Z=\emptyset$; indeed, if a
vertex of $Q$ was in $Z$ then it would be marked and so would its
neighbour in $Q$; thus we would have $Q=wx=w'x'$; but then one of $w',
x'$ would be in $Z\cap C(T)$, a contradiction.  Also note that if $w'=
w$ then this vertex is in $P\cap C(T)\cap N(Z)$ so, by
(\ref{eq:eqct}), $w'=w=u$.  On the other hand $x'=x$ is not possible
as $v\not\in Q$.  It follows that $Q$ has length at least $2$.

By the definition of marked vertices, there exists a subpath $Z'$ of
$Z$ with endvertices $z_i, z_j$ such that $z_iw, z_jx$ are edges (each
of $i<j$, $i=j$, $i>j$ is possible).  We choose $Z'$ minimal with that
property; so its internal vertices, if any, miss both $w, x$, and
consequently $H=Q\cup Z'$ is a hole of length at least $4$, and of
course it is an even hole.  Note that there is no vertex of $C(T)$ in
$H\setminus Q[w', x']$, by the definition of $w', x'$ and because
$Z\subseteq V\setminus(T\cup C(T))$.  Moreover, the $T$-edges in $H$
are exactly the $T$-edges in $Q$.

If $w', x'$ are not consecutive along $P$, then $H$ has two
non-adjacent vertices of $C(T)$ and yet it has an odd number of
$T$-edges, a contradiction to Lemma~\ref{lem:wh}.  So $w', x'$ must be
consecutive along $P$.

Put $k= \max \{i,j\}$ ($k\ge 1$).  Define a path $Y$ by setting $Y=
Z[z_{k+1}, z_n]\cup\{v\}$ if $k<n$ and $Y=\{v\}$ if $k=n$.  Note that
$H\cup Y$ is connected since $z_k$ is a vertex of $H$ adjacent to $Y$.
We claim that every vertex $z\in Y$ misses both $w', x'$.  Indeed, $v$
itself does miss both because $w',$ $x'$, $x$, $v$ are four distinct
vertices in that order along $P$; if $z\in Z[z_{k+1}, z_n]$ then $z$
misses $x'$ because $x'$ is not marked; moreover, if $z$ sees $w'$,
then $w'$ is marked, so $w'=w=u$, so $w'\in P\cap C(T)$, but then the
edge $zw'$ contradicts (\ref{eq:eqct}).  So the claim holds.  Now the
triple $(H, Y, T)$ contradicts Lemma~\ref{lem:dpg}, as there is no
vertex of $C(T)$ in $H\setminus \{w', x'\}$.  This completes the proof
that either $u'\in A$ or $v'\in B$ holds.

Finally, suppose that both $u'\in A$ and $v'\in B$ hold.  The path
$P'=P[u', v']$ has odd length and, since there is no edge between $A$
and $B$, this length is at least $3$.  Put $P'=u'u''\cdots v''v'$.
Applying the lemma to $P'$, we obtain that one of $u''\in A$ or
$v''\in B$ holds; but this contradicts the fact that $A$ and $B$ are
cliques and $P$ is chordless.  So exactly one of $u'\in A$ and $v'\in
B$ holds.  This completes the proof of the lemma.  \hfill$\Box$

We continue with the proof of Theorem~\ref{thm:main}.  Define a
relation $<_A$ on $A$ by setting $u<_A u'$ if and only if there exists
an odd chordless path from $u$ to a vertex of $B$ such that $u'$ is
the second vertex of that path.

\begin{lemma}\label{lem:antisym}
The relation $<_A$ is antisymmetric.
\end{lemma}
\emph{Proof.} Suppose on the contrary that there are vertices $u,v\in
A$ such that $u<_A v$ and $v<_A u$.  So there exists an odd chordless
path $P_u=u_0\cdots u_p$ such that $u=u_0$, $v=u_1$, $u_p\in B$, $p\ge
3$, $p$ odd, and there exists an odd chordless path $P_v=v_0\cdots
v_q$ such that $v=v_0$, $u=v_1$, $v_q\in B$, $q\ge 3$, $q$ odd.
Possibly $u_p=v_q$, and otherwise $u_p$ sees $v_q$ since $B$ is a
clique.  We observe that
\begin{eqnarray}\label{eq:nopupva}
\mbox{No vertex of $P_u\cup P_v\setminus\{u,v\}$ is in $A$,}
\end{eqnarray}
because such a vertex misses one of $u,v$ and $A$ is a clique.  Also,
\begin{eqnarray}\label{eq:nopupvb}
\mbox{No vertex of $P_u\cup P_v\setminus\{u_p,v_q\}$ is in $B$;}
\end{eqnarray}
indeed, since $B$ is a clique, a vertex of $(P_u\cup P_v\setminus
\{u_p,v_q\})\cap B$ can be only $u_{p-1}$ or $v_{q-1}$.  However, if
$u_{p-1}$ is in $B$ then $P_u[u_1, u_{p-1}]$ is an odd chordless path
from $A$ to $B$, and Lemma~\ref{lem:pq} implies either $u_2\in A$ (but
this contradicts (\ref{eq:nopupva})) or $u_{p-2}\in B$ (but this
contradicts that $B$ is a clique).  The case $v_{q-1}\in B$ is
similar.  So (\ref{eq:nopupvb}) is established.

Let $r$ be the smallest integer such that a vertex $u_r\in P_u
\setminus\{u_0, u_1\}$ has a neighbour in $P_v$, and let $s$ be the
smallest integer such that $u_rv_s$ is an edge, with $2\le s\le q$.
Such integers exist since $v_q$ itself has a neighbour in $P_u$.  Note
that $u_r\notin P_v$ and that the vertices $u_1, \ldots, u_r, v_1,
\ldots, v_s$ induce a hole $H$, so $r,s$ have the same parity, and
$u_r, v_s$ are different and adjacent.

Now we claim that we can assume that:
\begin{cadre}\label{cl:rs}
Either (a) $r=p$ and $s=q$, or (b) $r<p$, $s<q$ and $P_u[u_{r+1},
u_p]= P_v[v_{s+1}, v_q]$.
\end{cadre}
To prove this, let $t$ be the largest integer such that $u_rv_t$ is an
edge, with $2\le s\le t\le q$.

If $t-s$ is even then $P_u[u_1, u_r] \cup P_v[v_t, v_q]$ is an odd
chordless path from $A$ to $B$.  Its second vertex is $u_2$, and its
penultimate vertex $w$ is either $v_{q-1}$ (if $t< q$) or $u_r$ (if
$t=q$).  By Lemma~\ref{lem:pq} applied to that path, we have either
$u_2\in A$ (but this contradicts (\ref{eq:nopupva})) or $w\in B$.  By
(\ref{eq:nopupvb}) the latter is possible only if $w=u_r=u_p$ (so
$t=q$); in that case $P_v[v_1, v_s]\cup\{u_p\}$ is an odd chordless
path from $A$ to $B$, so, by Lemma~\ref{lem:pq}, we must have either
$v_2\in A$ (but this contradicts (\ref{eq:nopupva})) or $v_s\in B$;
this is possible only if $v_s=v_q$.  Thus we obtain case (a) of
(\ref{cl:rs}).

If $t-s$ is odd and $t\ge s+3$ then $P_v[v_1, v_s]\cup\{u_r\}\cup P_v
[v_t, v_q]$ is an odd chordless path from $A$ to $B$.  Its second
vertex is $v_2$, and its penultimate vertex $w$ is either $v_{q-1}$
(if $t< q$) or $u_r$ (if $t =q$).  By Lemma~\ref{lem:pq} applied to
that path, we must have either $v_2\in A$ (but this contradicts
(\ref{eq:nopupva})) or $w\in B$.  By (\ref{eq:nopupvb}) the latter is
possible only if $w=u_r=u_p$ (so $t=q$), thus $r, t$ are odd, but this
is impossible because $t-s$ is odd and $r-s$ is even.

The remaining possibility is when $t=s+1$.  Then the path $P_u[u_0,
u_r]\cup P_v[v_{s+1}, v_q]$ is odd and chordless, so it can play the
role of $P_u$, and we are in case (b).  So (\ref{cl:rs}) is proved.

Let $P$ be a chordless path defined as follows.  In case (a), set
$P=Z[z_2, z_n]$.  In case (b), set $P$ to be a chordless $(z_2,
u_{r+1})$-path contained in $Z[z_2, z_n] \cup P_u[u_{r+1}, u_p]$
(which induces a connected subgraph of $G$).  We observe that Lemma
\ref{lem:dpg} can be applied to the triple $(H, P, \{z_1\})$, indeed:
$H$ is a hole; $P$ is a chordless path; $\{z_1\}$ induces a
co-connected subgraph disjoint from $H\cup P$; $H\cup P$ induces a
connected subgraph; every vertex $z\in P$ misses both $u_1, v_1$ (when
$z\in Z[z_2, z_n]$ this is because $u_1, v_1\in A$, when $z\in P_u
[u_{r+1}, u_p]$ this is because $P_u, P_v$ are chordless paths and $r,
s\ge 2$); and $u_1, v_1, z_2$ are in $C(\{z_1\})$.  Now Lemma
\ref{lem:dpg} implies that some vertex of $H\setminus\{u_1, v_1\}$ is
in $C(\{z_1\})=N(z_1)$.  So $z_1$ has at least $3$ neighbours in $H$.
Call \emph{$z_1$-segment} any subpath of $H$, of length at least $1$,
whose endvertices are adjacent to $z_1$ and whose internal vertices
are not.  The conclusion of this paragraph is that {\it $H$ is
(edge-wise) partitioned into at least three $z_1$-segments.}

Now let us consider the vertices of $(H\setminus\{u_1, v_1\})\cap
C(T)$.  In case (a), both $u_r, v_s$ are in $C(T)$.  In case (b), we
can apply Lemma~\ref{lem:sgt} to the hole $H$, the path $P_u[u_{r+1},
u_p]$ and the set $T$, with respect to the edges $u_1v_1$ and
$u_rv_s$; Lemma~\ref{lem:sgt} implies that at least one of $u_r, v_s$
is in $C(T)$.  Thus, in either case (a) or (b), $H$ contains at least
three vertices of $C(T)$.  By Lemma~\ref{lem:wh}, we can conclude that
\emph{$H$ has an even number of $T$-edges.}

Observe that $u_1v_1$ is a $z_1$-segment that contains one $T$-edge.
It follows from this and the conclusion of the preceding two
paragraphs that there exists a $z_1$-segment $Q$ of $H$, different
from $u_1v_1$, that contains an odd number of $T$-edges.  Call $x, y$
the endvertices of $Q$, and call $x', y'$ respectively the first and
last vertex of $Q\cap C(T)$, so that $u_1, x, x', y', y, v_1$ lie in
this order along $H$.  Since $H$ has at least three $z_1$-segments, we
can assume that at least one of $x,y$, say $y$, is different from both
$u_1, v_1$; thus $y$ also misses one of $u_1, v_1$.

If $Q$ has length $1$, we have $Q=xy=x'y'$, so $y\in N(z_1)\cap C(T)$;
Lemma \ref{lem:abi} implies $y\in A\cup B\cup C(T\cup A\cup B)$; but
$y\in A\cup C(T\cup A\cup B)$ is impossible because $y$ misses one of
$u_1, v_1$, and $y\in B$ is impossible because $y\in N(z_1)$.  So $Q$
has length at least $2$.  Then $Q\cup \{z_1\}$ induces a hole $H_1$,
and Lemma~\ref{lem:wh} applied to the pair $H_1, T$ implies that
$Q[x', y']$ has length $1$.

If $x=u_1$, then $x'=u_1$ and $y'=u_2$, but then $\{v_1, z_1\}\cup
Q[u_2, y]$ induces an odd chordless path, with endvertices $v_1,
u_2\in C(T)$, that contains no $T$-edge, contradicting Lemma
\ref{lem:wp}.  So $x \not=u_1$.  Now consider the path induced by
$\{u_1, z_1\} \cup Q[y,y']$ and the path induced by $\{v_1, z_1\}\cup
Q[x, x']$; these paths are chordless (since $x\not=u_1$ and $y\not=
v_1$), their endvertices are in $C(T)$, they have no internal vertex
in $C(T)$, and one is them is odd (because $Q$ is even), so one of
them violates Lemma~\ref{lem:wp}.  This completes the proof of Lemma
\ref{lem:antisym}.  \hfill$\Box$

\begin{lemma}
The relation $<_A$ is transitive.
\end{lemma}
\emph{Proof.} Let $u, v, w$ be three vertices of $A$ such that $u <_A v
<_A w$.  Since $v <_A w$, there exists an odd chordless path $P=v_0
v_1 \cdots v_q$ with $v_0=v$, $v_1=w$, $v_q\in B$, $q$ odd, $q\ge 3$.

If $u$ has no neighbour along $P[v_2, v_q]$ then $\{u\}\cup P[v_1,
v_q]$ induces an odd chordless path to $B$, implying $u <_A w$ as
desired.  Now assume that $u$ has a neighbour $v_i$ along $P[v_2,
v_q]$, and let $i$ be the largest such integer ($2\le i\le q$).  We
have $i<q$ as there is no edge between $A$ and $B$.

If $i$ is odd ($3\le i\le q-2$), then $\{u\}\cup P[v_i, v_q]$ is an
odd chordless path with $u\in A$ and $v_q\in B$; applying Lemma
\ref{lem:pq} to this path, we have either $v_i\in A$ or $v_{q-1}\in
B$.  The former is impossible because $A$ is a clique; so $v_{q-1}\in
B$.  But then $\{v_0, u\}\cup P[v_i, v_{q-1}]$ induces an odd
chordless path to a vertex in $B$, which implies $v <_A u$,
contradicting Lemma~\ref{lem:antisym}.

If $i$ is even ($2\le i\le q-1$), then $\{v_0, u\}\cup P[v_i, v_q]$ is
an odd chordless path to a vertex in $B$, again implying $v <_A u$ and
contradicting Lemma~\ref{lem:antisym}.  This completes the proof of
the lemma.  \hfill$\Box$

{\bf Remark 3.} Determining the orders $<_A$ and $<_B$ and their
maximal elements can be done in polynomial time.  To do this, for any
three vertices $u,v\in A, w\in B$, look for a chordless $(v,w)$-path
in $G\setminus [(B \setminus\{w\})\cup (N(u)\setminus\{v\})]$.  If
such a path exists it must be even (by Lemma~\ref{lem:pq}) and its
existence implies $u<_A v$.  If no such path exists for any $w\in B$,
we have $u\not<_A v$ again by Lemma~\ref{lem:pq}.  For given $u,v\in
A, w\in B$, looking for such a path takes time $O(|E|)$, so, since
$A,B$ may both have size $O(|V|)$, the determination of $<_A$ and
$<_B$ takes time $O(|V|^3|E|)$.

\begin{lemma}\label{lem:absep}
Let $a$ be any maximal vertex of $(A, <_A)$ and $b$ be any maximal
vertex of $(B, <_B)$.  Then $\{a,b\}$ is a special even pair of $G$.
\end{lemma}
\emph{Proof.} Suppose that there exists an odd chordless $(a,b)$-path
$Q= a a' \cdots b'b$.  This path has length at least $3$ because there
is no edge between $A$ and $B$.  Lemma~\ref{lem:pq} implies either
$a'\in A$ (so $a<_A a'$) or $b'\in B$ (so $b<_B b'$), so in either
case the choice of $a$ or $b$ is contradicted.  So $\{a,b\}$ is an
even pair.

Now suppose that there exists a proper $(a,b)$-snake $S$ (with the
same notation as in the definition of snakes).  We may assume that
$S_1$ has length at least $1$, and we call $a_1$ the neighbour of $a'$
along $S_1$.  Since no vertex of $S$ sees both $a,b$, no vertex of $S$
is in $T$.  Note that $S_1\cup S_3\cup S_2$ and $S_1\cup S_4\cup S_2$
induce even paths since $\{a, b\}$ is an even pair.  Call $H$ the hole
induced by $S_3\cup S_4$.

We claim that:
\begin{eqnarray}
\mbox{No vertex of $S\setminus\{a,b\}$ is in $N(Z)\cap C(T)$.}
\label{eq:acdcd}
\end{eqnarray}
For suppose that there is a vertex $u$ of $N(Z) \cap C(T)$ in $S
\setminus\{a,b\}$.  By Lemma~\ref{lem:abi}, $u\in A\cup B\cup C(T \cup
A\cup B)$.  But $u \in C(T\cup A\cup B)$ is impossible because no
vertex of $S$ sees both $a,b$.  Therefore $u\in A\cup B$.  If $u\in
A$, $u$ must be the neighbour of $a$ along $S_1$ (since $A$ is a
clique); then $(S_1 \setminus\{a\})\cup S_3\cup S_2$ induces an odd
chordless path $P_u$ from $u\in A$ to $b$ (recall that $S_1\cup
S_3\cup S_2$ is an even chordless path); since the neighbour of $u$
along $P_u$ is not in $A$ (because $A$ is a clique),
Lemma~\ref{lem:pq} implies that the neighbour of $b$ along $P_u$ is in
$B$; but this means that $b$ is not maximal in $(B, <_B)$, a
contradiction.  If $u\in B$, $u$ must be either the neighbour of $b$
along $S_2$ (if $S_2$ has length at least $1$) or one of $c',d'$ (if
$S_2$ has length $0$), but in either case, an argument similar to the
case when $u\in A$ holds.  Thus (\ref{eq:acdcd}) is proved.

By Lemma~\ref{lem:sgp} applied to $S, T$, we know that:
\begin{eqnarray}
\mbox{At least two of $a',c,d$ and two of $b', c', d'$ are in
$C(T)$.}\label{eq:ald}
\end{eqnarray}

Since $a\neq a'$ (but $b=b'$ is possible), Facts (\ref{eq:acdcd}) and
(\ref{eq:ald}) imply:
\begin{cadre}
If one of $a',c,d$ is in $N(Z)$, then it is in $N(Z)\setminus C(T)$
and the other two are in $C(T)\setminus N(Z)$.  \newline If one of
$c', d'$ is in $N(Z)$, then it is in $N(Z)\setminus C(T)$ and the
other is in $C(T)\setminus N(Z)$, and $b'\in C(T)$.\label{cd:ctnz}
\end{cadre}
Moreover, $a'\not\in Z$ for otherwise one of $c,d$ would be in
$N(Z)\cap C(T)$.

Now we define a path $P$ as follows.  Let $a''$ be the vertex of
$N(Z)$ closest to $a_1$ along $S_1\setminus \{a'\}$, and let $b''$ be
the vertex of $N(Z)$ closest to $b'$ along $S_2$ (vertices $a'', b''$
exist because of $a,b$).  Pick $z_i \in Z\cap N(a'')$ and $z_j\in Z
\cap N(b'')$ such that the path $Z[z_i, z_j]$ is as short as possible
(each of $i<j$, $i=j$, $i>j$ is possible).  Put $P=S_1[a_1, a'']\cup
Z[z_i, z_j]\cup S_2[b'', b']$; so $P$ is a chordless $(a_1, b')$-path.
By Lemma~\ref{lem:dpg} applied to the triple $(H, P, \{a'\})$, some
vertex $z$ of $P$ sees one of $c, d$, and we can pick $z$ closest to
$a_1$ along $P$ and assume up to symmetry that $z$ sees $c$.  The
definition of $S$ implies $z \not\in S_1[a_1, a'']\cup S_2$, so $z\in
Z[z_i, z_j]$.  By (\ref{cd:ctnz}), we have $c\in N(z) \setminus C(T)$
and $a', d\in C(T) \setminus N(Z)$.  Thus $a'$ has no neighbour along
$P\setminus \{a_1\}$, for such a neighbour could only be in $Z[z_i,
z_j]$ (by the definition of $S$), and then we would have $a'\in N(Z)$,
a contradiction.  In other words,
\begin{eqnarray}
\mbox{$P\cup\{a'\}$ is a chordless path.}\label{eq:pa1}
\end{eqnarray}
By Lemma~\ref{lem:dpg} applied to the triple $(H, P[z, b']\setminus
\{b'\}, \{b'\})$, some vertex $y$ of the path $P[z, b'] \setminus
\{b'\}$ sees one of $c', d'$.  By the definition of $S$ and $P$, we
have $y\in P[z, z_j]=Z[z, z_j]$.  Thus, and by (\ref{cd:ctnz}), we
know that: %
\begin{cadre}
Exactly one of $c',d'$ is in $C(T)\setminus N(Z)$, the other is in
$N(Z)\setminus C(T)$, and $b'\in C(T)$.\label{cd:cpdpt}
\end{cadre}

We note that:
\begin{eqnarray}
\mbox{$cc'\not\in E$.}\label{eq:ccnz}
\end{eqnarray}
For suppose $cc'\in E$.  If $c'\in C(T)$, then $a'cc'$ is an outer
path of length two, a contradiction.  If $c'\not\in C(T)$, then, by
(\ref{cd:cpdpt}), $a'cc'b'$ is an odd outer path, a contradiction.  So
(\ref{eq:ccnz}) holds.

Call $H_1$ the cycle induced by $P[a_1, z]\cup\{a', c\}$, which has
length at least $4$.  By (\ref{eq:pa1}) and by the choice of $z$,
$H_1$ is a hole.  Call $S'_4$ the path $S_4\cup\{c'\}$.  We claim
that:
\begin{eqnarray}
\mbox{There is no edge between $P[a_1, z]$ and $S'_4$.}\label{eq:nowx}
\end{eqnarray}
For suppose that there is an edge $xw$ with $x\in P[a_1, z]$ and $w\in
S'_4$.  We have $w\not=d$ since $d\notin N(Z)$ and $d$ has no
neighbour on $S_1\setminus\{a'\}$.  No vertex of $S'_4\setminus \{d\}$
is adjacent to $a'$ or $c$ by the definition of $S$ and
by~(\ref{eq:ccnz}).  But then the triple $(H_1, S'_4 \setminus\{d\},
\{d\})$ violates Lemma~\ref{lem:dpg}.  So (\ref{eq:nowx}) holds.

We know that some vertex $z'$ of $P[z, z_j]$ has a neighbour $\delta$
in $S'_4$ (because of $y$), so we can pick $z'$ closest to $z$ along
$P[z, z_j] = Z[z, z_j]$ and pick $\delta\in S'_4\cap N(z')$ closest to
$d$ along $S'_4$.  Thus $P[z, z']\cup S'_4[\delta, d]$ is a chordless
path.  We have $\delta\not=d$ since $z'\in Z$ and $d\notin N(Z)$.  Let
us consider the cycle $H_2$ induced by $P[a_1, z']\cup S'_4[\delta,
d]\cup\{a'\}$, which has length at least $4$.  Suppose that $H_2$ has
a chord.  The definition of $S, P, z', \delta$ and the fact that
$P\cup\{a'\}$ and $P[z, z']\cup S'_4 [\delta, d]$ are chordless imply
that the only possible chords in $H_2$ are of the type $wx$ with $w\in
S_4$ and $x\in Z[z_i, z]$.  But this is forbidden by (\ref{eq:nowx}).
So $H_2$ is an even hole.  Now, along $H_2$, vertex $c$ has three
neighbours $a', d, z$; hence, by Lemma~\ref{lem:wh}, $H_2$ has an even
number of $c$-edges (edges whose two endvertices are neighbours of
$c$).  One of these is $a'd$.  Obviously there is no $c$-edge along
$S'_4[d, \delta]$.  Also $\delta z'$ cannot be a $c$-edge, for that
would imply $\delta=c'$ and $cc'\in E$, contradicting (\ref{eq:ccnz}).
Thus all the $c$-edges of $H_2$ different from $a'd$ (and there is an
odd number of these) lie in $P[z, z']$.  Call $z''$ the neighbour of
$c$ closest to $z'$ along $P[z, z']$, so that all the $c$-edges of
$H_2$ different from $a'd$ lie in $P[z, z'']$.  By Lemma~\ref{lem:wp}
applied to $P[z, z'']$ and $\{c\}$, we obtain that:
\begin{eqnarray}
\mbox{$P[z, z'']$ has odd length.}
\end{eqnarray}

By Lemma~\ref{lem:wh}, the number of $T$-edges in $H_1$ is either $1$
or even.  Suppose it is $1$; so the vertices of $H_1\cap C(T)$ are
$a'$ and $a_1$ since $c\notin C(T)$.  By Lemma~\ref{lem:dpg} applied
to the triple $(H_1, P[b', z]\setminus\{z\}, T)$, some vertex of
$P[b', z]\setminus\{z\}$ sees one of $a', a_1$; but this contradicts
(\ref{eq:pa1}).  So:
\begin{eqnarray}
\mbox{$H_1$ has an even number of $T$-edges.}
\end{eqnarray}

\begin{figure}[htb]
\unitlength=0.4cm
\begin{center}
\begin{picture}(14,12)
\multiput(4,0)(8,0){2}{\circle*{0.3}}
\multiput(4,0)(8,0){2}{\circle{0.7}}
\multiput(2,2)(12,0){2}{\circle{0.7}}
\multiput(4,4)(2,0){5}{\circle*{0.3}}
\multiput(0,2)(2,0){2}{\circle*{0.3}}
\multiput(3,8)(1,0.6){5}{\circle*{0.3}}
\multiput(5.4,0)(1.4,0){3}{\circle*{0.3}}
\put(14,2){\circle*{0.3}}
\put(9.6,0){\circle*{0.3}}
\put(11.2,0){\circle*{0.3}}
\put(8.8,10.4){\circle*{0.3}}
\put(-1,1.5){$a_1$}
\put(1.5,2.5){$a'$}
\put(4.5,3){$c$}
\put(11.3,3){$c'$}
\put(4.5,0.5){$d$}
\put(11.3,0.5){$d'$}
\put(10,0.5){$\delta'$}
\put(8.4,0.5){$\delta$}
\put(14.5,2.5){$b'$}
\put(2.2,8.2){$z$}
\put(5.3,10){$z''$}
\put(7,10.8){$z'$}
\put(8.8,10.8){$y$}
\put(13,8){$P$}
\multiput(4,0)(0,4){2}{\line(1,0){8}}
\multiput(4,0)(8,0){2}{\line(0,1){4}}
\multiput(2,2)(10,-2){2}{\line(1,1){2}}
\multiput(2,2)(10,2){2}{\line(1,-1){2}}
\put(0,2){\line(1,0){2}}
\put(3,8){\line(5,3){4}}
\put(4,4){\line(-1,4){1}}
\put(4,4){\line(0,1){4.6}}
\put(4,4){\line(1,3){1.9}}
\put(7,10.4){\line(2,-5){4.2}}
\put(7,10.4){\line(1,-4){2.6}}
\put(7,10.4){\line(1,0){1.8}}
\put(12,4){\line(-1,2){3.2}}
\qbezier(0,2)(0,6)(3,8)
\qbezier(8.8,10.4)(14,8.8)(14,2)
\end{picture}
\end{center}
\vspace{0.6cm} \caption{Lemma~\protect{\ref{lem:absep}} in case $d'\in
C(T)$.  Circled vertices are in $C(T)$.}\label{fig:final}
\end{figure}
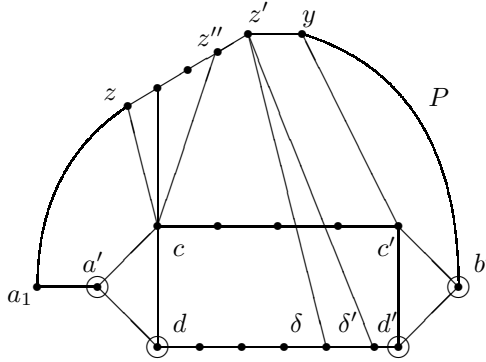

In view of (\ref{cd:cpdpt}), let us call $\gamma'$ whichever of $c',
d'$ is in $C(T)$ and define a path $S''_4$ as follows.  If $\delta=c'$
(so $\gamma'=d'$), we put $S''_4= c'd'$.  If $\delta\in S_4$ then let
$\delta'$ be the neighbour of $z'$ along $S_4$ that is closest to
$d'$, and put $S''_4 =S'_4[\delta', \gamma']$.  Now consider the path
$Q_1$ induced by $\{a', c\}\cup P[z'', z']\cup S''_4$ and the path
$Q_2$ induced by $\{a'\}\cup P[a_1, z'] \cup S''_4$.  These paths are
chordless by (\ref{eq:pa1}), (\ref{eq:ccnz}), (\ref{eq:nowx}) and the
definition of $z', z'', \delta'$.  Their endvertices are $a'$ and
$\gamma'$, which are both in $C(T)$.  Since $H_1$ is an even hole and
$P[z, z'']$ has odd length, the paths $Q_1, Q_2$ have lengths of
different parity.  The numbers of $T$-edges in $Q_1$ and in $Q_2$ have
the same parity, because $H_1$ has an even number of $T$-edges and
$cz$ is not a $T$-edge (as $c\notin C(T)$), and because $P[z, z'']$
contains no $T$-edge (as $P[z, z'']\subseteq Z$).  Thus one of $Q_1,
Q_2$ violates Lemma~\ref{lem:wp}, a contradiction.  \hfill$\Box$

The conclusion of Lemma~\ref{lem:absep} completes the proof of Theorem
\ref{thm:main}.  \hfill$\Box$

\section{Some consequences}\label{sec:more}
\subsubsection*{Complexity and Optimal Colorings}

Along the proof of Theorem~\ref{thm:main} we made three remarks on the
polynomial complexity of finding various sets, paths or orders related
to that proof.  In total these remarks show that the proof of
Theorem~\ref{thm:main} is a polynomial-time algorithm which, given any
graph $G=(V,E)$ that is not a clique, either returns a special even
pair of $G$ or answers that $G$ contains an odd hole, an antihole of
length at least $5$ or a prism.  We note that Case~1 of the proof
leads to an iterative call of the algorithm on a subgraph of $G$; this
may happen $O(|V|$) time.  On the other hand Case~2 does not lead to
iterating and therefore happens only once during the execution of the
whole algorithm.  Thus, a rough estimate of the complexity of this
algorithm is $O(|V|^3|E|)$.

As suggested in the Introduction, we can use this algorithm to color
every graph $G$ in ${\cal A}$ using no more than $\omega(G)$ colors.
If $G$ is not a disjoint union of cliques, we use the above algorithm
to get a special even pair $\{x,y\}$ of $G$, then iterate the
procedure with the graph $G/xy$.  If $G$ is a disjoint union of
cliques an $\omega(G)$-coloring can be produced trivially.  As there
are at most $O(|V|)$ iterations, this coloring algorithm has time
complexity $O(|V|^4|E|)$.

\subsubsection*{Subclasses of ${\cal A}$}

The class ${\cal A}$ contains several families of perfect graphs for
which the existence of an even pair was already proved, in a specific
way for each such family
\cite{hayhoamaf90,her90,herdew88,mey87,rus00}; we will not recall the
definition of all these families here, see the survey \cite{epsbook}.
Let us however make a few remarks about two such families.

1.  A graph $G$ is \emph{perfectly orderable} if it admits a
\emph{perfect ordering}, that is, an ordering of its vertices such
that, for every induced subgraph $G'$ of $G$, using the greedy
coloring method on the vertices of $G'$ along the induced ordering
produces an optimal coloring of $G'$.  Chv\'atal \cite{chv84}
introduced perfectly orderable graphs and proved that they are
perfect.  Hertz and de Werra \cite{herdew88} showed that every
perfectly orderable graph $G$ is perfectly contractile by proving that
if $G$ is not a clique it has an even pair $\{x,y\}$ such that $G/xy$
is also perfectly orderable.  Their proof assumes that a perfect
ordering for $G$ is given.  However, determining if a graph is
perfectly orderable is an NP-complete problem \cite{midpfe90}, so
there is probably no efficient way to find an even pair in a perfectly
orderable graph using that method if no perfect ordering is given.
Our result bypasses this difficulty.  A drawback is that if $G$ is
perfectly orderable, not a clique, and $\{a,b\}$ is the even pair
produced by our algorithm, we cannot certify that $G/ab$ is also
perfectly orderable, only that it is in ${\cal A}$.

2.  A graph $G$ is \emph{weakly triangulated} if $G$ and
$\overline{G}$ contain no hole of length at least $5$.  Hayward
\cite{hay85} proved that weakly triangulated graphs are perfect, and
later Hayward, Ho\`ang and Maffray \cite{hayhoamaf90} proved that they
are perfectly contractile using the following definition and theorem.
A \emph{ $2$-pair} is a pair of vertices $\{u,v\}$ such that all
chordless $(u,v)$-paths have length $2$.
\begin{theorem}[\cite{hayhoamaf90}]\label{thm:wt}
If $G$ is a weakly triangulated graph and not a clique, then:
\begin{enumerate}
\item
$G$ has a $2$-pair;
\item
For any $2$-pair $\{a,b\}$ of $G$, $G/ab$ is weakly triangulated.
\end{enumerate}
\end{theorem}
We show here an alternate proof of the first item of Theorem
\ref{thm:wt}.  If $G$ is a disjoint union of cliques then any two
non-adjacent vertices form a $2$-pair.  If $G$ is not a disjoint union
of cliques, we can mimick the proof of Theorem~\ref{thm:main}:
consider any maximal interesting set $T$; by Lemma~\ref{lem:wt}, there
is no outer path with respect to $T$, so we are in Case~1; in Case~1,
we can assume that the induction hypothesis provides a $2$-pair of
$G[C(T)]$; then the same arguments as in Case~1 imply that it is also
a $2$-pair of $G$.  In other words, when the input graph $G$ is weakly
triangulated, our algorithm produces a $2$-pair of $G$.  \hfill$\Box$

\section*{Acknowledgement}

We are grateful to the two referees and to Celina de Figueiredo and
Cl\'audia Villela Maciel for their careful reading and comments on the
paper.  We thank all of them for pointing out a mistake in the last
lemma in an earlier version of the manuscript.

Part of this work was done while the authors were guests of the
American Institute of Mathematics (AIM, http://aimath.org) during the
Perfect Graph Workshop organized by Paul Seymour and Robin Thomas.  We
thank AIM and the organizers for the invitation and generous
support.


%

\vfill
\begin{flushright}April 23 2004\end{flushright}

\end{document}